\documentclass[10pt,reqno]{amsart}

\usepackage{amsmath,amssymb}
\usepackage[dvips]{graphics}
\usepackage{epsfig}

\newtheorem{thm}{Theorem}[section]

\newtheorem{cor}[thm]{Corollary}

\newtheorem{exa}[thm]{Example}
\newtheorem{defi}[thm]{Definition}

\numberwithin{equation}{section}

\newtheorem{rek}[thm]{Remark}
\theoremstyle{definition}

\newcommand\ben{\begin{enumerate}}
\newcommand\een{\end{enumerate}}
\newcommand\bi{\begin{itemize}}
\newcommand\ei{\end{itemize}}

\newcommand{\twocase}[5]{#1 \begin{cases} #2 & \text{#3}\\ #4
&\text{#5} \end{cases}   }

\newcommand\be{\begin{equation}}
\newcommand\ee{\end{equation}}
\newcommand\bea{\begin{eqnarray}}
\newcommand\eea{\end{eqnarray}}

\newcommand{\R}{\mathbb{R}}

\newcommand{\Z}{\mathbb{Z}}

\newcommand{\ga}{\alpha}     
\newcommand{\gb}{\beta}      
\newcommand{\gep}{\epsilon}  

\renewcommand{\mod}{\;\operatorname{mod}}

\newcommand{\E}{{\mathbb E}} 

\numberwithin{equation}{section}

\begin{document}

\title{Order Statistics and Benford's Law}
\author{Steven J. Miller}
\email{Steven.J.Miller@williams.edu}\address{Department of Mathematics and Statistics, Williams College, Williamstown, MA 01267}\author{Mark
J. Nigrini} \email{nigrini@tcnj.edu}
\thanks{
We thank Ted Hill, Christoph Leuenberger, Daniel Stone and the referees for numerous helpful comments.
The first author was
partially supported by NSF grant DMS-0600848.}
\address{Accounting and Information Systems, School of Business, The College of New Jersey, Ewing, NJ 08628} 
\date{\today}



\begin{abstract} Fix a base $B>1$ and let $\zeta$ have the standard exponential
distribution; the distribution of digits of $\zeta$ base $B$ is
known to be very close to Benford's Law. If there exists a $C$ such
that the distribution of digits of $C$ times the elements of some set is the same as that of $\zeta$, we say that set exhibits
shifted exponential behavior base $B$ (with a shift of $\log_B C
\bmod 1$). Let $X_1, \dots, X_N$ be independent identically
distributed random variables. If the $X_i$'s are drawn from the
uniform distribution on $[0,L]$, then as $N\to\infty$ the distribution
of the digits of the differences between adjacent order statistics
converges to shifted exponential behavior (with a shift of
$\log_B L/N \bmod 1$).
By differentiating the cumulative distribution
function of the logarithms modulo 1, applying Poisson Summation and
then integrating the resulting expression, we derive rapidly
converging explicit formulas measuring the deviations from Benford's
Law.
Fix a $\delta \in (0,1)$ and choose $N$ independent random variables
from any compactly supported distribution with uniformly bounded first and second derivatives and a second order Taylor series expansion at each point. The distribution of digits of any
$N^\delta$ consecutive differences \emph{and} all $N-1$ normalized
differences of the order statistics exhibit shifted exponential behavior.
We derive conditions on the probability density which determine
whether or not the distribution of the digits of all the
un-normalized differences converges to Benford's Law, shifted exponential behavior, or oscillates between the two, and show that the Pareto distribution leads to oscillating behavior.

\end{abstract}

\maketitle


\section{Introduction}

Benford's Law gives the expected frequencies of the digits in many
tabulated data. It was first observed by Newcomb in the 1880s, who
noticed that pages of numbers starting with a $1$ in logarithm tables
were significantly more worn than those starting with a $9$. In 1938
Benford \cite{Ben} observed the same digit bias in a variety of
phenomenon. From his observations he postulated that in many data
sets more numbers began with a 1 than with a 9; his investigations (with 20,229 observations) supported his belief.
See \cite{Hi1,Rai} for a description and history and
\cite{Hu} for an extensive bibliography.

For any base $B>1$ we may uniquely write a positive $x\in\R$ as $x =
M_B(x)\cdot B^k$, where $k\in \Z$ and $M_B(x)$ (called the mantissa)
is in $[1,B)$. A sequence of positive numbers $\{a_n\}$ is
\textbf{Benford base $B$} if the probability of observing a mantissa
of $a_n$ base $B$ of at most $s$ is $\log_B s$. More precisely, for $s \in [1,B]$ we have \be
\lim_{N \to \infty} \frac{ \#\{n \le N: \text{$1 \le M_B(a_n) \le
s$} \} }{N} \ =\ \log_B s. \ee Benford behavior for continuous
functions\footnote{If the functions are not positive, we study the distribution of the digits of the absolute value of the function.} are defined analogously. Thus base $10$ the probability of
observing a first digit of $d$ is $\log_{10} (d+1) - \log_{10} (d)$,
implying that about $30\%$ of the time the first digit is a $1$.

We can prove many mathematical systems follow Benford's law, ranging from recurrence relations \cite{BrDu} to $n!$ \cite{Dia} to iterates of power, exponential and
rational maps and Newton's method \cite{Hi2,BBH,BH} to chains of random variables and hierarchical Bayesian models \cite{JKKKM} to values of $L$-functions near the critical line to
characteristic polynomials of random matrix ensembles and
iterates of the $3x+1$-Map \cite{KonMi,LS} to products of random variables \cite{MN}; we also see Benford's law in a variety of natural systems, such as atomic physics \cite{Pa}, biology \cite{CLTF} and geology \cite{NM1}. Applications of Benford's Law range from
rounding errors in computer calculations (see page 255 of
\cite{Knu}) to detecting tax (see \cite{Nig1,Nig2}) and voter fraud (see \cite{Me}).


This work is motivated by two observations (see Remark \ref{rek:motivation} for more details). First, since Benford's
seminal paper, many investigations have shown that amalgamating data
from different sources leads to Benford behavior; second, many
standard probability distributions are close to Benford behavior.
We investigate the distribution of digits of differences of adjacent
ordered random variables. For any $\delta < 1$, if we study at most $N^\delta$ consecutive differences of a data set of size $N$, the resulting
distribution of leading digits depends very weakly on the underlying
distribution of the data, and closely approximates Benford's Law. We
then investigate whether or not studying all the differences lead to
Benford behavior; this question is inspired by the first observation
above, and has led to new tests for data integrity (see \cite{NM2}). These tests are quick and easy to apply, and have successfully detected problems with some data sets, thus providing a practical application of our main results.

To prove our results requires analyzing the distribution of digits
of independent random variables drawn from the standard exponential,
and quantifying how close the distribution of digits of a random
variable with the standard exponential distribution is to Benford's
Law. Leemis, Schmeiser and Evans \cite{LSE} have observed that the
standard exponential is quite close to Benford's Law; this was
proved by Engel and Leuenberger \cite{EL}, who showed that the
maximum difference in the cumulative distribution function from
Benford's Law (base $10$) is at least .029 and at most .03. We
provide an alternate proof of this result in the appendix using a different technique, as well as showing that there is no base $B$ such that the standard exponential distribution is Benford base $B$ (Corollary \ref{cor:neverbenford}).

Both proofs apply Fourier analysis to periodic functions. In
\cite{EL} the main step (their equation (5)) is interchanging an
integration and a limit. Our proof is based on applying Poisson
Summation to the derivative of the cumulative distribution function
of the logarithms modulo $1$, $F_B$. Benford's Law is equivalent to
$F_B(b) = b$, which by calculus is the same as $F_B'(b) = 1$ and
$F_B(0)=0$. Thus studying the deviation of $F_B'(b)$ from $1$ is a
natural way to investigate the deviations from Benford behavior. We hope the details of these calculations may be of use to
others in investigating related problems (Poisson Summation
has been fruitfully used by Kontorovich-Miller \cite{KonMi} and Jang-Kang-Kruckman-Kudo-Miller \cite{JKKKM} in
proving many systems are Benford; see also
\cite{Pin}).

\subsection{Definitions}

A sequence $\{a_n\}_{n=1}^\infty \subset [0,1]$ is equidistributed
if \be \lim_{N\to \infty} \frac{ \# \{n : n \leq N,\ a_n\in [a,b]\}
}{N} \ = \ b-a\ee for all $[a,b]\subset [0,1]$. Similarly a
continuous random variable on $[0,\infty)$ whose probability density
function is $p$ is equidistributed modulo $1$ if \be
\lim_{T\to\infty} \frac{\int_0^T \chi_{a,b}(x)p(x)dx}{\int_0^T
p(x)dx} \ = \ b-a \ee for any $[a,b] \subset [0,1]$, where
$\chi_{a,b}(x) = 1$ for $x\bmod 1 \in [a,b]$ and $0$ otherwise.

A positive sequence (or values of a function) is Benford base $B$ if and only
if its base $B$ logarithms are equidistributed modulo $1$; this
equivalence is at the heart of many investigations of Benford's Law;
see \cite{Dia,MT-B} for a proof.

We use the following notation for the various error terms:

\ben

\item Let $\mathcal{E}(x)$ denote an
error of at most $x$ in absolute value; thus $f(b) = g(b) +
\mathcal{E}(x)$ means $|f(b) - g(b)| \le x$.

\item big-Oh notation: For $g(x)$ a non-negative function, we say
$f(x) = O(g(x))$ if there exists an $x_0$ and a $C>0$ such that, for
all $x \ge x_0$, $|f(x)| \le C g(x)$.

\een

The following theorem is the starting point for investigating the
distribution of digits of order statistics.

\begin{thm}\label{thm:main} Let $\zeta$ have
the standard (unit) exponential distribution: \be\label{eq:stexpdistz} {\rm
Prob}\left(\zeta \in [\alpha,\beta]\right) \ = \ \int_\alpha^\beta
e^{-t}dt, \ \ \ \ [\alpha,\beta] \in [0,\infty). \ee For $b \in
[0,1]$, let $F_B(b)$ be the cumulative distribution function of
$\log_B\zeta \bmod 1$; thus $F_B(b) := {\rm Prob}(\log_B\zeta \bmod
1 \in [0,b])$. Then for all $M \ge 2$ \bea\label{eq:thmeqmainFprime}
F_B'(b) &\ =\ & 1 + 2 \sum_{m=1}^\infty {\rm Re}\left(e^{-2\pi i
mb} \Gamma\left(1+\frac{2\pi i m}{\log B}\right)\right) \nonumber\\
& = & 1 + 2\sum_{m=1}^{M-1} {\rm Re}\left(e^{-2\pi i mb}
\Gamma\left(1+\frac{2\pi i m}{\log B}\right)\right)\nonumber\\ & & \
\ \ \ \ + \ \mathcal{E}\left(4\sqrt{2}\pi c_1(B)
e^{-(\pi^2-c_2(B))M/\log B}\right), \eea where $c_1(B), c_2(B)$ are
constants such that for all $m \ge M \ge 2$ we have \bea e^{2\pi^2
m/\log B} -
e^{-2\pi^2 m/\log B} & \ \ge \ & e^{2\pi^2 m / \log B} / c_1^2(B) \nonumber\\
m/\log B & \le & e^{2c_2(B)m/\log B} \nonumber\\ 1 -
e^{-(\pi^2-c_2(B))M/\log B} & \ge & 1/\sqrt{2}. \eea For $B \in
[e,10]$ we may take $c_1(B) =\sqrt{2}$ and $c_2(B) = 1/5$, which
give \bea\label{eq:qerwerw} {\rm Prob}(\log\zeta \bmod 1 \in [a,b])&
\ = \ & b - a\ +\ \frac{2r}{\pi}\cdot \sin(\pi (b+a) + \theta)\cdot
\sin( \pi(b-a))\nonumber\\ & & \ \ \ \ \  + \ \mathcal{E}\left(6.32
\cdot 10^{-7}\right), \eea with $r \approx 0.000324986$, $\theta
\approx 1.32427186$, and \bea\label{eq:qerwerwten}& & {\rm
Prob}(\log_{10}\zeta \bmod 1 \in [a,b]) \ = \  b - a\ +\
\frac{2r_1}{\pi} \sin(\pi (b+a) - \theta_1)\cdot \sin(
\pi(b-a))\nonumber\\ & &\ \ \ \ \ \ \ \ \ \ \ \ \ \ \ \ \ \ \ \ -
\frac{r_2}{\pi} \sin(2\pi (b+a) + \theta_2)\cdot \sin( 2\pi(b-a))\ +
\ \mathcal{E}(8.5 \cdot 10^{-5}),\eea with \bea r_1 & \ \approx \ &
0.0569573, \ \ \ \ \theta_1 \ \approx \ 0.8055888 \nonumber\\ r_2 &
\ \approx \ & 0.0011080, \ \ \ \ \theta_2 \ \approx \ 0.1384410.
\eea
\end{thm}

The above theorem was proved in \cite{EL}; we provide an alternate
proof in the appendix. As remarked earlier, our technique consists of applying Poisson Summation to the derivative of the cumulative distribution function of the logarithms modulo $1$; it is then very natural and easy to compare deviations from the resulting distribution and the uniform distribution (if a data set satisfies Benford's law, then the distribution of its logarithms is uniform). Our series expansions are obtained by applying properties of the Gamma function.

\begin{defi}[Exponential Behavior, Shifted Exponential Behavior] Let $\zeta$ have the standard exponential distribution,
and fix a base $B$. If the distribution of the digits of a set is
the same as the distribution of the digits of $\zeta$, then we say
the set exhibits exponential behavior (base $B$). If there is
a constant $C>0$ such that the distribution of digits of all
elements multiplied by $C$ is exponential behavior, then we say
the system exhibits shifted exponential behavior (with shift of
$\log_B C \bmod 1$).
\end{defi}

We briefly describe the reasons behind this notation. One important
property of Benford's Law is that it is invariant under rescaling;
many authors have used this property to characterize Benford
behavior. Thus if a data set is Benford base $B$ and we fix a
positive number $C$, so is the data set obtained by multiplying each
element by $C$. This is clear if, instead of looking at the
distribution of the digits, we study the distribution of the base
$B$ logarithms modulo $1$. Benford's Law is equivalent to the
logarithms modulo $1$ being uniformly distributed (see for instance \cite{Dia,MT-B}); the effect of
multiplying all entries by a fixed constant simply translates the
uniform distribution modulo $1$, which is again the uniform
distribution.

The situation is different for exponential behavior. Multiplying
all elements by a fixed constant $C$ (where $C \neq B^k$ for some
$k\in\Z$) does not preserve exponential behavior; however, the
effect is easy to describe. Again looking at the logarithms, exponential behavior is equivalent to the base $B$ logarithms modulo $1$
having a specific distribution which is almost equal to the uniform
distribution (at least if the base $B$ is not too large). Multiplying by a fixed constant $C \neq B^k$ shifts
the logarithm distribution by $\log_B C \mod 1$.

\subsection{Results for Differences of Orders Statistics}

We consider a simple case first, and show how the more general case
follows. Let $X_{1}, \dots, X_{N}$ be independent identically
distributed from the uniform distribution on $[0,L]$. We consider
$L$ fixed and study the limit as $N\to\infty$. Let $X_{1:N}, \dots,
X_{N:N}$ be the $X_i$'s in increasing order. The $X_{i:N}$ are
called the order statistics, and satisfy $0 \le X_{1:N} \le X_{2:N} \le
\cdots \le X_{N:N} \le L$. We investigate the distribution of the
leading digits of the differences between adjacent $X_{i:N}$'s,
$X_{i+1:N} - X_{i:N}$. For convenience we periodically continue the
data and set $X_{i+N:N} = X_{i:N}+L$. As we have $N$ differences in
an interval of size $L$, on average $X_{i+1:N} - X_{i:N}$ is of size
$L/N$, and it is sometimes easier to study the normalized
differences \be\label{eq:normalizeddiffzxs}
Z_{i;N} \ = \ \frac{X_{i+1:N}-X_{i:N}}{L/N}.
\ee As the $X_i$'s are drawn from a uniform distribution, it is a
standard result that as $N\to\infty$ the $Z_{i;N}$'s are independent
random variables, each having the standard exponential distribution.
Thus as $N\to\infty$ the probability that $Z_{i;N} \in [a,b]$ tends to $\int_a^b e^{-t}dt$. See \cite{DN,Re} for proofs.

For uniformly distributed random variables, if we know the
distribution of $\log_B Z_{i;N} \bmod 1$ then we can immediately determine the
distribution of the digits of the $X_{i+1:N} - X_{i:N}$ base $B$
because \be\label{eq:erere} \log_B Z_{i;N} \ = \
\log_B\left(\frac{X_{i+1:N}-X_{i:N}}{L/N}\right)\ = \ \log_B
(X_{i+1:N} - X_{i:N}) - \log_B (L/N). \ee

As the $Z_{i;N}$ are independent with the standard exponential
distribution as $N\to\infty$ if the $X_i$ are independent uniformly
distributed, the behavior of the digits of the differences
$X_{i+1:N} - X_{i:N}$ is an immediate consequence of Theorem
\ref{thm:main}:

\begin{thm}[Shifted Exponential Behavior of Differences of Independent
Uniformly Distributed Random Variables]\label{thm:unifdistdiffab}
Let $X_{1}, \dots, X_{N}$ be independently distributed from the
uniform distribution on $[0,L]$, and let $X_{1:N}, \dots, X_{N:N}$
be the $X_i$'s in increasing order. As $N\to\infty$ the distribution
of the digits (base $B$) of the differences $X_{i+1:N} - X_{i:N}$
converges to shifted exponential behavior, with a shift of
$\log_B (L/N) \bmod 1$.
\end{thm}

A similar result holds for other distributions.

\begin{thm}[Shifted Exponential Behavior of Subsets of Differences of Independent Random
Variables] \label{thm:gennicedistrdiff} Let $X_{1},\dots, X_{N}$ be
independent, identically distributed random variables whose density $f(x)$ has
a second order Taylor series at each point with first and second
derivatives uniformly bounded, and let the $X_{i:N}$'s be the
$X_i$'s in increasing order. Fix a $\delta \in (0,1)$. Then as
$N\to\infty$ the distribution of the digits (base $B$) of $N^\delta$
consecutive differences $X_{i+1:N} - X_{i:N}$ converges to shifted
exponential behavior, provided the $X_{i:N}$'s are from a region
where $f(x)$ is non-zero.
\end{thm}

The key ingredient in this generalization is that the techniques
which show that the differences between uniformly distributed random
variables become independent exponentially distributed random
variables can be modified to handle more general distributions.

We restricted ourselves to a subset of all consecutive spacings
because the normalization factor changes throughout the domain. The
shift in the shifted exponential behavior depends on which set of
$N^\delta$ differences we study, coming from the variations in the
normalizing factors. Within a bin of $N^\delta$ differences the
normalization factor is basically constant, and we may approximate
our density with a uniform distribution. It is possible for these
variations to cancel and yield Benford behavior for the digits of
all the un-normalized differences. Such a result is consistent with
the belief that amalgamation of data from many different
distributions becomes Benford; however, this is not always the case
(see Remark \ref{rek:paretodangerabb}). From Theorem \ref{thm:main}
and Theorem \ref{thm:gennicedistrdiff} we obtain

\begin{thm}[Benford Behavior for all the Differences of Independent Random
Variables]\label{thm:shiftsab} Let $X_{1},\dots, X_{N}$ be
independent, identically distributed random variables whose density $f(x)$ is compactly supported and has
a second order Taylor series at each point with first and second
derivatives uniformly bounded. Let the $X_{i:N}$'s be the
$X_i$'s in increasing order, $F(x)$ be the cumulative distribution function for $f(x)$, and 
fix a $\delta \in (0,1)$. Let $I(\gep,\delta,N) = [\gep
N^{1-\delta}, N^{1-\delta} - \gep N^{1-\delta}]$. For each fixed
$\gep \in (0, 1/2)$, assume that
\begin{itemize}
\item $f(F^{-1}(kN^{\delta-1})$ is not too small for $k\in I(\gep,\delta,N)$:
\be\label{eq:firstassumpnotsmall} \lim_{N\to\infty}\ \max_{k \in
I(\gep,\delta,N)}\ \frac{\min(N^{-(\gep + \delta/2)},
N^{\delta-1})}{f(F^{-1}(kN^{\delta-1}))} \ = \ 0; \ee
\item $\log_B f(F^{-1}(kN^{\delta-1}) \bmod 1$ is
equidistributed: for all $[\alpha,\beta] \subset
[0,1]$\be\label{eq:thmwhatneedequidistrabb} \lim_{N\to\infty} \frac{
\#\{k \in I(\gep,\delta,N): \log_B f(F^{-1}(kN^{\delta-1})) \bmod 1
\in [\alpha, \beta]\}}{N^{\delta}} \ = \ \beta - \alpha. \ee
 \end{itemize}
Then if $\gep > \max(0, 1/3-\delta/2)$ and $\gep < \delta/2$, the distribution of the digits of the $N-1$ differences
$X_{i+1:N} - X_{i:N}$ converges to Benford's Law (base $B$) as
$N\to\infty$.
\end{thm}

\begin{rek}\label{rek:paretodangerabb}
The conditions of Theorem \ref{thm:shiftsab} are usually not
satisfied. We are unaware of any situation where \eqref{eq:thmwhatneedequidistrabb} holds; we have included Theorem \ref{thm:shiftsab} to give a sufficient condition of what is required to have Benford's law satisfied \emph{exactly}, and not just approximately. In Lemma \ref{lem:paretonoshiftgood} we show the
conditions fail for the Pareto distribution, and the limiting
behavior oscillates between Benford and a sum of shifted exponential behavior.\footnote{If several data sets each exhibit shifted
exponential behavior but with distinct shifts, then the amalgamated data set is closer to Benford's Law than any of the original data sets.
This is apparent by studying the logarithms modulo $1$. The differences between these densities and Benford's law will look like the plot on the right in Figure \ref{fig:Pareto500000binAllof100Scaled} (except, of course, that different shifts will result in shifting the plot modulo $1$). The key observation is that the unequal shifts mean we do not have reinforcements from the peaks of the modulo $1$ densities being aligned, and thus the amalgamation will decrease the maximum deviations.} The
arguments generalize to many densities whose cumulative
distribution functions have tractable closed-form expressions (for
example, exponential, Weibull, or $f(x) = e^{-e^x} e^x$).
\end{rek}

The situation is very different if instead we study normalized
differences \be\label{eq:normalizeddiffzxstilde} \widetilde{Z}_{i:N}
\ = \ \frac{X_{i+1:N} - X_{i:N}}{1 / N f(X_{i:N})}; \ee note if
$f(x) = 1/L$ is the uniform distribution on $[0,L]$,
\eqref{eq:normalizeddiffzxstilde} reduces to
\eqref{eq:normalizeddiffzxs}.

\begin{thm}[Shifted Exponential Behavior for
All the Normalized Differences of Independent Random
Variables]\label{thm:shiftaddreinforce} Assume the probability
distribution $f$ satisfies the conditions of Theorem
\ref{thm:shiftsab} and \eqref{eq:firstassumpnotsmall} and
$\widetilde{Z}_{i;N}$ is as in \eqref{eq:normalizeddiffzxstilde}.
Then as $N\to\infty$ the distribution of the digits of the
$\widetilde{Z}_{i:N}$ converges to shifted exponential behavior.
\end{thm}

\begin{rek} Appropriately scaled, the distribution of the digits of the
differences is universal, and is the exponential behavior of
Theorem \ref{thm:main}. Thus Theorem \ref{thm:shiftaddreinforce}
implies that the natural quantity to study is the normalized
differences of the order statistics, not the differences. See also
Remark \ref{rek:universalabnormdiff}. With additional work we could study densities with unbounded support and show that, through truncation, we can get arbitrarily close to shifted exponential behavior. \end{rek}

\begin{rek}\label{rek:motivation} The main motivation for this work is the need for improved ways
of assessing the authenticity and integrity of scientific and
corporate data. Benford's Law has been successfully applied to
detecting income tax, corporate and voter fraud (see \cite{Me,Nig1,Nig2}); in
\cite{NM2} we use these results to derive new statistical tests to
examine data authenticity and integrity.
Early applications of these tests to financial data showed that it could detect errors in data
downloads, rounded data, and inaccurate ordering of data.  These attributes are not easily observable from
an analysis of descriptive statistics, and detecting these errors can help managers avoid costly decisions based on erroneous data.\end{rek}

The paper is organized as follows. We prove Theorem \ref{thm:main}
in Appendix \ref{sec:approachviaderivs} by using Poisson summation
to analyze $F_B'(b)$. Theorem \ref{thm:unifdistdiffab} follows
from results for the order statistics of independent uniform
variables; the proof of Theorem \ref{thm:gennicedistrdiff} is
similar, and given in \S\ref{subsec:orderniceprobdistr}. In
\S\ref{sec:distrdigdiffallos} we prove Theorems \ref{thm:shiftsab}
and \ref{thm:shiftaddreinforce}.



\section{Proofs of Theorem \ref{thm:unifdistdiffab} and
\ref{thm:gennicedistrdiff}}\label{subsec:orderniceprobdistr}

Theorem \ref{thm:unifdistdiffab} is a consequence of the fact that
the normalized differences between the order statistics drawn from
the uniform distribution converge to being independent standard
exponentials. The proof of Theorem \ref{thm:gennicedistrdiff}
proceeds similarly. Specifically, over a short enough region any distribution with a second order Taylor series at each point with first and second derivatives uniformly bounded is well-approximated by a uniform distribution.

To prove Theorem \ref{thm:gennicedistrdiff}, it suffices to show
that if $X_{1}, \dots, X_{N}$ are drawn from a sufficiently nice
distribution, then for any fixed $\delta \in (0,1)$ the limiting
behavior of the order statistics of $N^\delta$ adjacent $X_i$'s
becomes Poissonian (i.e., the $N^\delta-1$ normalized differences
converge to being independently distributed from the standard
exponential). We prove this below for compactly supported distributions $f(x)$ that have
a second order Taylor series at each point with the first and second
derivatives uniformly bounded, and when the $N^\delta$ adjacent
$X_i$'s are from a region where $f(x)$ is bounded away from zero.

For each $N$, consider intervals $[a_N,b_N]$ such that
$\int_{a_N}^{b_N} f(x)dx = N^\delta / N$; thus the proportion of the
total mass in such intervals is $N^{\delta-1}$. We fix such an
interval for our arguments. For each $i \in \{1,\dots, N\}$ let \be\label{eq:defnwi}
\twocase{w_i \ = \ }{1}{if $X_i \in [a_N,b_N]$}{0}{otherwise.} \ee
Note $w_i$ is $1$ with probability $N^{\delta-1}$ and $0$ with
probability $1-N^{\delta-1}$; $w_i$ is a binary indicator random
variable, telling us whether or not $X_i \in [a_N,b_N]$. Thus
\bea\label{eq:tretete} \E\left[\sum_{i=1}^N w_i\right] \ =\
N^\delta, \ \ \ \ \  {\rm Var}\left(\sum_{i=1}^N w_i\right)  \ = \
N^{\delta} \cdot (1 - N^{\delta-1}). \eea Let $M_N$ be the number of
$X_i$ in $[a_N,b_N]$, and let $\beta_N$ be \emph{any} non-decreasing sequence tending to infinity (in the course of the proof, we will find we may take any sequence with $\beta_N = o(N^{\delta/2})$). By \eqref{eq:tretete} and the Central Limit
Theorem (which we may use as the $w_i$'s satisfy the Lyapunov condition), with probability tending to $1$ we have \be\label{eq:MNCLT}
M_N \ = \ N^\delta + O(\beta_N N^{\delta/2}).\ee

We assume that in the interval $[a_N,b_N]$ there exist constants $c$
and $C$ such that whenever $x\in [a_N,b_N]$,  $0 < c < f(x) < C <
\infty$; we assume these constants hold for all regions investigated and for all $N$.\footnote{If our distribution has unbounded support, for any $\gep > 0$ we can truncate it on both sides so that the omitted probability is at most $\gep$. Our result is then trivially modified to being within $\gep$ of shifted exponential behavior.} Thus \be c \cdot (b_N-a_N) \ \le \ \int_{a_N}^{b_N} f(x)dx
\ = \ N^{\delta - 1} \ \le \ C (b_N-a_N), \ee implying that
$b_N-a_N$ is of size $N^{\delta - 1}$. If we assume $f(x)$ has at
least a second order Taylor expansion, then
\bea\label{eq:taylorexpf} f(x) & \ = \ & f(a_N) + f'(a_N) (x-a_N) +
O((x-a_N)^2)\nonumber\\ & \ = \ & f(a_N) + f'(a_N) (x-a_N) +
O(N^{2\delta-2}).\eea \emph{As we are assuming the first and second
derivatives are uniformly bounded, as well as $f$ being bounded away
from zero in the intervals under consideration, all big-Oh constants
below are independent of $N$.} Thus \be\label{eq:diffbnan} b_N -
a_N\ =\ \frac{N^{\delta - 1}}{f(a_N)} + O(N^{2\delta-2}). \ee

We now investigate the order statistics of the $M_N$ of the $X_i$'s
that lie in $[a_N,b_N]$. We know $\int_{a_N}^{b_N} f(x)dx =
N^{\delta - 1}$; by setting $g_N(x) = f(x) N^{1-\delta}$ then
$g_N(x)$ is the conditional density function for $X_i$, given that $X_i \in [a_N, b_N]$. Thus $g_N(x)$ integrates to 1, and for $x \in [a_N,b_N]$ we have
\be\label{eq:taylorexpansiongnx} g_N(x) \ = \ f(a_N)\cdot N^{1-\delta} + f'(a_N)(x-a_N)\cdot
N^{1-\delta} + O(N^{\delta-1}).\ee

We have an interval of size $N^{\delta-1}/f(a_N) +
O(N^{2\delta-2})$, and $M_N = N^\delta + O(\beta_N N^{\delta/2})$ of the
$X_i$ lying in the interval (remember the $\beta_N$ are any non-decreasing sequence tending to infinity). Thus with probability tending to 1, the average spacing between
adjacent ordered $X_i$ is \bea\label{eq:anbnavespaceingos}
\frac{N^{\delta - 1}/f(a_N) + O(N^{2\delta-2})}{M_N}& \ = \ &
(f(a_N)N)^{-1} + N^{-1} \cdot O(\beta_N N^{-\delta/2} + N^{\delta-1}); \eea in particular, we see we must choose $\beta_N = o(N^{\delta/2})$.  As $\delta \in (0,1)$, if we fix a $k$ such that $X_k \in [a_N,b_N]$
then we expect the next $X_i$ to the right of $X_k$ to be about
$\frac{t}{Nf(a_N)}$ units away, where $t$ is of size $1$. For a
given $X_k$ we can compute the conditional probability that the next $X_i$ is between
$\frac{t}{Nf(a_N)}$ and $\frac{t+\Delta t}{Nf(a_N)}$ units to the
right: it is simply the difference of the probability that all the other $M_N-1$ of the 
$X_i$'s in $[a_N,b_N]$ are not in the interval $[X_k, X_k +
\frac{t}{Nf(a_N)}]$ and the probability that all other $X_i$ in
$[a_N,b_N]$ are not in the interval $[X_k, X_k + \frac{t+\Delta
t}{Nf(a_N)}]$; note we are using the wrapped interval $[a_N,b_N]$.

Some care is required in these calculations. We have a conditional probability as we are assuming both $X_k \in [a_N, b_N]$ and that exactly $M_N$ of the $X_i$ are in $[a_N, b_N]$. Thus these probabilities depend on two random variables, namely $X_k$ and $M_N$. This is not a problem in practice, however (for example, $M_N$ is tightly concentrated about its mean value).

Recalling our expansion for $g_N(x)$ (and that $b_N-a_N =
N^{\delta-1}/f(a_N) + O(N^{2\delta-2})$ and $t$ is of size $1$), 
after simple algebra we find that, with probability tending to 1, for a given $X_k$ and $M_N$ the first probability is \bea \left(1 -
\int_{X_k}^{X_k + \frac{t}{Nf(a_N)}} g_N(x)dx\right)^{M_N-1}. \eea 
The above integral equals $tN^\delta + O(N^{-1})$ (use the Taylor series expansion in \eqref{eq:taylorexpansiongnx} and note that the interval $[a_N, b_N]$ is of size $O(N^{\delta-1})$).
Using \eqref{eq:MNCLT}, is easy to see that this is a.s. equal to
\bea 
\left(1 - \frac{t + O(N^{\delta-1} +
\beta_N N^{-\delta/2})}{M_N}\right)^{M_N-1}. \eea  We therefore find that as $N \to \infty$ the probability that $M_N-1$ of the $X_i$'s ($i \neq k$) are in $[a_N, b_N] \setminus [X_k, X_k + t/Nf(a_N)]$, conditioned on $X_k$ and $M_N$, converges to $e^{-t}$.\footnote{Some care is required, as the exceptional set in our a.s. statement can depend on $t$. This can be surmounted by taking expectations with respect to our conditional probabilities and applying the dominated convergence theorem.}

The calculation of the second probability, the conditional probability that the $M_N - 1$ other $X_i$'s in
$[a_N,b_N]$ are not in the interval $[X_k, X_k + \frac{t+\Delta
t}{Nf(a_N)}]$, given $X_k$ and $M_N$, follows analogously by replacing $t$ with $t+\Delta t$ in the previous argument. We thus find that this probability is $e^{-(t+\Delta t)}$. As \be \int_t^{t+\Delta t} e^{-u}du \ = \ e^{-t} - e^{-(t+\Delta t)},\ee we find that the density of the difference between adjacent order statistics tends to the standard (unit) exponential density; thus the proof of Theorem \ref{thm:gennicedistrdiff} now follows from Theorem \ref{thm:unifdistdiffab}.

%


\section{Proofs of Theorems \ref{thm:shiftsab} and
\ref{thm:shiftaddreinforce}}\label{sec:distrdigdiffallos}

We generalize the notation from \S\ref{subsec:orderniceprobdistr}.
Let $f(x)$ be any distribution with a second order Taylor series at each point with first and second derivatives uniformly bounded, and let
$X_{1:N}, \dots, X_{N:N}$ be the order statistics. We fix a $\delta
\in (0,1)$, and for $k \in \{1,\dots,N^{1-\delta}\}$ we consider bins $[a_{k;N}, b_{k;N}]$ such that \be
\int_{a_{k;N}}^{b_{k;N}} f(x) dx \ = \ N^\delta / N \ = \
N^{\delta-1}; \ee there are $N^{1-\delta}$ such bins. By the Central
Limit Theorem (see \eqref{eq:MNCLT}), if $M_{k;N}$ is the number of
order statistics in $[a_{k;N}, b_{k;N}]$ then provided that $\gep > \max(0, 1/3 - \delta/2)$ with probability
tending to $1$ we have \be\label{eq:mknsimultepsilon} M_{k;N} \ = \ N^{\delta}
+ O(N^{\gep+\delta/2});\ee of course, we also require $\gep < \delta/2$, as otherwise the error term is larger than the main term.

\begin{rek} Before we considered just one fixed
interval; as we are studying $N^{1-\delta}$ intervals
simultaneously, we need the $\gep$ in the exponent so that with high
probability all intervals have to first order $N^\delta$ order
statistics. For the arguments below, it would have sufficed to have
an error of size $O(N^{\delta-\gep})$. We thank the referee for pointing out that $\gep > 1/3 - \delta/2$, and provide his argument in Appendix \ref{sec:analyzingintervalssimult}.  \end{rek}

Similar to \eqref{eq:anbnavespaceingos}, the average spacing between
adjacent order statistics in $[a_{k;N}, b_{k;N}]$ is
\be\label{eq:erereqazaq} (f(a_{k;N})N)^{-1} + N^{-1} \cdot
O(N^{-(\gep+\delta/2)} + N^{\delta-1}). \ee Note
\eqref{eq:erereqazaq} is the generalization of \eqref{eq:erere}; if
$f$ is the uniform distribution on $[0,L]$ then $f(a_{k;N}) = 1/L$.
By Theorem \ref{thm:gennicedistrdiff}, as $N\to\infty$ the
distribution of digits of the differences in each bin converges to
shifted exponential behavior; however, the variation in the
average spacing between bins leads to bin-dependent shifts in the
shifted exponential behavior.

Similar to \eqref{eq:erere}, we can study the distribution of digits
of the differences of the normalized order statistics. If $X_{i:N}$
and $X_{i+1:N}$ are in $[a_{k;N}, b_{k;N}]$ then
\bea\label{eq:zinscaledeff} Z_{i;N} & = & (X_{i+1:N} - X_{i:N})
\Big/ \left((f(a_{k;N})N)^{-1} + N^{-1} \cdot O(N^{-(\gep+\delta/2)}
+ N^{\delta-1})\right) \nonumber\\ \log_B Z_{i;N} &=&  \log_B
(X_{i+1:N} - X_{i:N}) + \log_B N  -\log_B \left(f(a_{k;N})^{-1} +
O(N^{-(\gep+\delta/2)} + N^{\delta-1})\right). \nonumber\\ \eea Note
we are using the \emph{same} normalization factor for all
differences between adjacent order statistics in a bin. Later we
show we may replace $f(a_{k;N})$ with $f(X_{i:N})$. As we study all
$X_{i+1:N} - X_{i:N}$ in the bin $[a_{k;N}, b_{k;N}]$, it is useful
to rewrite the above as \bea\label{eq:explogxi+1xi} \log_B(X_{i+1:N}
- X_{i:N}) & = & \log_B Z_{i;N} - \log_B N + \log_B
\left(f(a_{k;N})^{-1} + O(N^{-(\gep+\delta/2)} +
N^{\delta-1})\right). \nonumber\\ \eea We have $N^{1-\delta}$ bins,
so $k \in \{1,\dots,N^{1-\delta}\}$. As we only care about the
limiting behavior, we may safely ignore the first and last bins. We
may therefore assume each $a_{k;N}$ is finite, and $a_{k+1;N} =
b_{k;N}$.\footnote{Of course, we know both quantities are finite as we assumed our distribution has compact support. We remove the last bins to simplify generalizations to non-compactly supported distributions.}

Let $F(x)$ be the cumulative distribution function for $f(x)$. Then
\be F(a_{k;N})\ =\ (k-1) N^\delta / N \ = \ (k-1)N^{\delta-1}.\ee
For notational convenience we relabel the bins so that $k \in
\{0,\dots, N^{1-\delta}-1\}$; thus $F(a_{k;N}) = k N^{\delta-1}$.

We now prove our theorems which determine when these bin-dependent
shifts cancel (yielding Benford behavior), or reinforce (yielding
sums of shifted exponential behavior).

\begin{proof}[Proof of Theorem \ref{thm:shiftsab}]
There are approximately $N^{\delta}$ differences in each bin
$[a_{k;N}, b_{k;N}]$. By Theorem \ref{thm:gennicedistrdiff}, the
distribution of the digits of the differences in each bin converges
to shifted exponential behavior. As we assume the first and
second derivatives of $f$ are uniformly bounded, the big-Oh
constants in \S\ref{subsec:orderniceprobdistr} are independent of
the bins. The shift in the shifted exponential behavior in each
bin is controlled by the last two terms on the right hand side of
\eqref{eq:explogxi+1xi}. The $\log_B N$ shifts the shifted exponential behavior in each bin equally. The bin-dependent shift is
controlled by the final term, \bea\label{eq:bindependentshiftterms}
& & \log_B \left(f(a_{k;N})^{-1} + O(N^{-(\gep+\delta/2)} +
N^{\delta-1})\right)\nonumber\\ &\ = \ & - \log_B f(a_{k;N}) +
\log_B\left(1 + \frac{\min(N^{-(\gep+\delta/2)},
N^{\delta-1})}{f(a_{k;N})}\right). \eea

Thus each of the $N^{1-\delta}$ bins exhibits shifted exponential
behavior, with a bin-dependent shift composed of the two terms in
\eqref{eq:bindependentshiftterms}. By
\eqref{eq:firstassumpnotsmall}, the $f(a_{k;N})$ are not small
compared to $\min(N^{-(\gep+\delta/2)}, N^{\delta-1})$, and hence
the second term $\log_B\left(1 + \frac{\min(N^{-(\gep+\delta/2)},
N^{\delta-1})}{f(a_{k;N})}\right)$ is negligible. In particular,
this factor depends only very weakly on the bin, and tends to zero
as $N\to\infty$.

Thus the bin-dependent shift in the shifted exponential
behavior is approximately $-\log_B f(a_{k;N})$ $=$ $-\log_B
f(F^{-1}(kN^{\delta-1}))$. If these shifts are equidistributed
modulo $1$, then the deviations from Benford behavior cancel, and
the shifted exponential behavior of each bin becomes Benford
behavior for \emph{all} the differences.
\end{proof}

\begin{rek} Consider the case when the density is a uniform distribution on some
interval. Then all $f(F^{-1}(kN^{\delta-1}))$ are equal, and each
bin has the same shift in its shifted exponential behavior. These
shifts therefore reinforce each other, and the distribution of all
the differences is also shifted exponential behavior, with the same shift.
This is observed in numerical experiments; see Theorem
\ref{thm:unifdistdiffab} for an alternate proof.
\end{rek}

We analyze the assumptions of Theorem \ref{thm:shiftsab}. The
condition from \eqref{eq:firstassumpnotsmall} is easy to check, and
is often satisfied. For example, if the probability density is a
finite union of monotonic pieces and is zero only finitely often,
then \eqref{eq:firstassumpnotsmall} holds. This is because for $k
\in I(\gep,\delta,N)$, $F^{-1}(kN^{\delta-1}) \in
[F^{-1}(\gep),F^{-1}(1-\gep)]$ and is therefore independent of $N$
(if $f$ vanishes finitely often, we need to remove small
sub-intervals from $I(\gep,\delta,N)$, but the analysis proceeds
similarly). The only difficulty is basically a probability
distribution with intervals of zero probability. Thus
\eqref{eq:firstassumpnotsmall} is a mild assumption.

If we choose any distribution \emph{other than} a uniform
distribution, then $f(x)$ is not constant; however,
\eqref{eq:thmwhatneedequidistrabb} need not hold (i.e., $\log_B
f(a_{k;N}) \bmod 1$ need not be equidistributed as $N\to\infty$).
For example, consider a Pareto distribution with minimum value $1$
and exponent $a>0$. The density is \be \twocase{f(x) \ = \ }{a
x^{-a-1}}{if $x \ge 1$}{0}{otherwise.}\ee The Pareto distribution is
known to be useful in modeling natural phenomena, and for
appropriate choices of exponents yields approximately Benford
behavior (see \cite{NM1}).

\begin{exa}\label{lem:paretonoshiftgood}
If $f$ is a Pareto distribution with minimum value $1$ and exponent
$a>0$, then $f$ does not satisfy the second condition of Theorem
\ref{thm:shiftsab}, equation \eqref{eq:thmwhatneedequidistrabb}.

To see this, note that the
cumulative distribution function of $f$ is $F(x) = 1 - x^{-a}$. As
we only care about the limiting behavior, we need only study $k \in
I(\gep,\delta,N)= [\gep N^{1-\delta}, N^{1-\delta} - \gep
N^{1-\delta}]$. Therefore $F(a_{k;N}) = kN^{\delta-1}$ implies that
\be a_{k;N} \ = \ (1 - kN^{\delta-1})^{-1/a}, \ \ \ f(a_{k;N}) \ = \
a (1 - kN^{\delta-1})^{\frac{a+1}{a}}. \ee 

The condition from \eqref{eq:firstassumpnotsmall} is satisfied,
namely \be \lim_{N\to\infty} \max_{k\in I(\gep,\delta,N)}
\frac{\min(N^{-(\gep+\delta/2)}, N^{\delta-1})}{f(a_{k;N})} \ = \
\lim_{N\to\infty} \max_{k\in
I(\gep,\delta,N)}\frac{\min(N^{-(\gep+\delta/2)},
N^{\delta-1})}{a(kN^{\delta-1})^{(a+1)/a}} \ = \ 0, \ee as $k$ is of
size $N^{1-\delta}$.

Let $j = N^{1-\delta} - k \in I(\gep,\delta,N)$. Then the
bin-dependent shifts are \bea \log_B f(a_{k;N}) &\ = \ &
\frac{a+1}{a}\log_B(1 - kN^{\delta-1}) + \log_B a \nonumber\\ & \ =
\ & \frac{a+1}{a} \log_B (j N^{1-\delta} ) + \log_B a \nonumber\\
&=& \log_B \left(j^{(a+1)/a}\right) + \log_B\left(a
N^{(1-\delta)(a+1)/a}\right). \eea Thus, for a Pareto distribution
with exponent $a$, the distribution of \emph{all} the differences
becomes Benford if and only if $j^{(a+1)/a}$ is Benford. This
follows from the fact that a sequence is Benford if and only if its
logarithms are equidistributed. For fixed $m$, $j^{m}$ is \emph{not}
Benford (see for example \cite{Dia}), and thus the condition from
\eqref{eq:thmwhatneedequidistrabb} fails.
\end{exa}

\begin{rek} We chose to study a Pareto distribution because the
distribution of digits of a random variable drawn from a Pareto
distribution converges to Benford behavior (base 10) as $a\to 1$;
however, the digits of the differences do not tend to Benford (or
shifted exponential) behavior. A similar analysis holds for many
distributions with good closed-form expressions for the cumulative
distribution function. In particular, if $f$ is the density of an
exponential or Weibull distribution (or $f(x) = e^{-e^x}e^x$), then
$f$ does not satisfy the second condition of Theorem
\ref{thm:shiftsab}, equation \eqref{eq:thmwhatneedequidistrabb}.
\end{rek}

Modifying the proof of Theorem \ref{thm:shiftsab} yields our result
on the distribution of digits of the normalized differences.

\begin{proof}[Proof of Theorem \ref{thm:shiftaddreinforce}]
If $f$ is the uniform distribution, there is nothing to prove. For
general $f$, rescaling the differences eliminates the bin-dependent
shifts. Let \be\label{eq:widetildezin} \widetilde{Z}_{i:N} \ = \
\frac{X_{i+1:N} - X_{i:N}}{1 / N f(X_{i:N})}. \ee In Theorem
\ref{thm:shiftsab} we use the same scale factor for all differences
in a bin; see \eqref{eq:zinscaledeff}. As we assume the first and
second derivatives of $f$ are uniformly bounded,
\eqref{eq:taylorexpf} and \eqref{eq:diffbnan} imply that for
$X_{i:N} \in [a_{k;N}, b_{k;N}]$, \bea f(X_{i:N}) &\ = \ &
f(a_{k;N}) + O\left(b_{k;N}-a_{k;N}\right) \nonumber\\ & = &
f(a_{k;N}) + O\left(\frac{N^{\delta - 1}}{f(a_{k;N})} +
N^{2\delta-2}\right), \eea and the big-Oh constants are independent
of $k$. As we assume $f$ satisfies \eqref{eq:firstassumpnotsmall},
the error term is negligible.

Thus our assumptions on $f$ imply that $f$ is basically constant on
each bin, and we may replace the local rescaling factor $f(X_{i:N})$
with the bin rescaling factor $f(a_{k;N})$. Thus each bin of
normalized differences has \emph{the same} shift in its shifted
exponential behavior. Therefore all the shifts reinforce, and the
digits of all the normalized differences exhibit shifted exponential behavior as $N\to\infty$.
\end{proof}

\begin{figure}[ht]
\begin{center}
\scalebox{.78}{\includegraphics{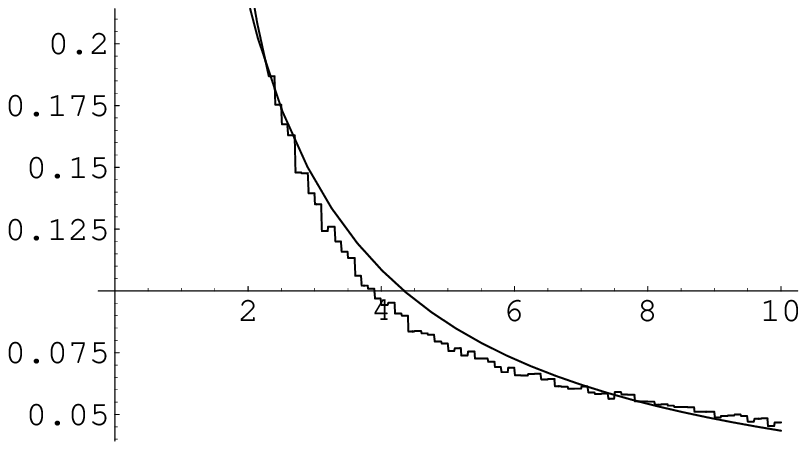}}\
\
\scalebox{.78}{\includegraphics{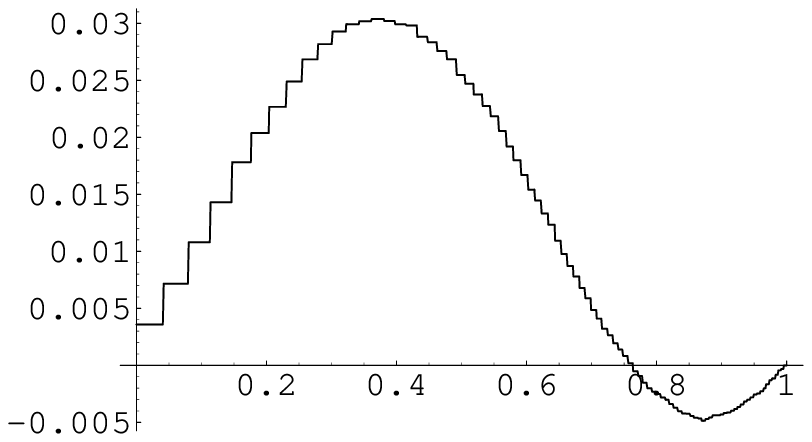}}
\caption{\label{fig:Pareto500000binAllof100Scaled}All 499,999
differences of adjacent order statistics from 500,000 independent
random variables from the Pareto distribution with minimum value and
variance $1$. (left) Observed digits of scaled differences of adjacent
random variables versus Benford's law; (right) Scaled observed minus
Benford's Law (cumulative distribution of base $10$ logarithms). }
\end{center}\end{figure}

As an example of Theorem \ref{thm:shiftaddreinforce}, in Figure
\ref{fig:Pareto500000binAllof100Scaled} we consider 500,000
independent random variables drawn from the Pareto distribution with
exponent \be a \ = \
\frac{4+\sqrt[3]{19-3\sqrt{33}}+\sqrt[3]{19+3\sqrt{33}}}3 \ee (we
chose $a$ to make the variance equal $1$). We study the distribution
of the digits of the differences in base 10. The amplitude is about
$.018$, which is the amplitude of the shifted exponential
behavior of Theorem \ref{thm:main} (see the equation in Theorem 2 of
\cite{EL} or \eqref{eq:thmeqmainFprime} of Theorem \ref{thm:main}).

\begin{rek}\label{rek:universalabnormdiff} The universal behavior of Theorem
\ref{thm:shiftaddreinforce} suggests that if we are interested in
the behavior of the digits of all the differences, the natural
quantity to study is the \emph{normalized} differences. For any distribution with uniformly bounded first and second derivatives and a second order Taylor series expansion at each point, we obtain shifted exponential
behavior.
\end{rek}


\appendix


\section{Proof of Theorem \ref{thm:main}}\label{sec:approachviaderivs}

To prove Theorem \ref{thm:main} it suffices to study the
distribution of $\log_B \zeta \bmod 1$ when $\zeta$ has the standard
exponential distribution; see \eqref{eq:stexpdistz}. We have the
following useful chain of equalities. Let $[a,b] \subset [0,1]$.
Then \bea\label{eq:startkeyequalities} {\rm Prob}(\log_B \zeta \bmod
1 \in [a,b]) & \ = \ & \sum_{k=-\infty}^\infty {\rm Prob}(\log_B
\zeta \in [a+k, b+k]) \nonumber\\ & = & \sum_{k=-\infty}^\infty {\rm
Prob}(\zeta \in [B^{a+k},B^{b+k}])\nonumber\\  & = &
\sum_{k=-\infty}^\infty \left(e^{-B^{a+k}} -
e^{-B^{b+k}}\right).\eea It suffices to investigate
\eqref{eq:startkeyequalities} in the special case when $a=0$, as the
probability of any interval $[\ga,\gb]$ can always be found by
subtracting the probability of $[0,\ga]$ from $[0,\gb]$. We are
therefore led to studying, for $b\in [0,1]$, the cumulative
distribution function of $\log_B \zeta \bmod 1$:
\bea\label{eq:problogzisum1} F_B(b) \ := \ {\rm Prob}(\log_B \zeta
\bmod 1 \in [0,b]) & \ = \ & \sum_{k=-\infty}^\infty
\left(e^{-B^{k}} - e^{-B^{b+k}}\right).\eea This series expansion
converges rapidly, and Benford behavior for $\zeta$ is equivalent to
the rapidly converging series in \eqref{eq:problogzisum1} equalling
$b$ for all $b$.

As Benford behavior is equivalent to $F_B(b)$ equals $b$ for all
$b\in [0,1]$, it is natural to compare $F_B'(b)$ to $1$. If the
derivative were identically $1$ then $F_B(b)$ would equal $b$ plus
some constant. However, \eqref{eq:problogzisum1} is zero when $b=0$,
which implies that this constant would be zero. It is hard to
analyze the infinite sum for $F_B(b)$ directly. By studying the
derivative $F_B'(b)$ we find a function with an easier Fourier
transform than the Fourier transform of $e^{-B^u} - e^{-B^{b+u}}$,
which we then analyze by applying Poisson Summation.

We use the fact that the derivative of the infinite sum $F_B(b)$ is
the sum of the derivatives of the individual summands. This is
justified by the rapid decay of the summands; see, for example,
Corollary 7.3 of \cite{La}. We find
\bea\label{eq:FprimeBderivfirsttime} 
F_B'(b) \ = \ \sum_{k=-\infty}^\infty e^{-B^{b+k}} B^{b+k} \log B \
= \ \sum_{k=-\infty}^\infty e^{-\beta B^k} \beta B^k \log B, \eea
where for $b \in [0,1]$ we set $\beta  =  B^b$.

Let $H(t) = e^{-\beta B^t} \beta B^t \log B$; note $\beta \ge 1$. As
$H(t)$ is of rapid decay in $t$, we may apply Poisson Summation (see
for example \cite{SS}). Thus \be \sum_{k=-\infty}^\infty H(k) \ = \
\sum_{k=-\infty}^\infty \widehat{H}(k), \ee where $\widehat{H}$ is
the Fourier Transform of $H$: $\widehat{H}(u)  =
\int_{-\infty}^\infty H(t) e^{-2\pi i tu} dt$. Therefore \bea
F_B'(b) \ =\  \sum_{k=-\infty}^\infty H(k)\ =\
\sum_{k=-\infty}^\infty \widehat{H}(k)\ =\ \sum_{k=-\infty}^\infty
\int_{-\infty}^\infty e^{-\beta B^t} \beta B^t \log B \cdot e^{-2\pi
i tk} dt. \eea Let us change variables by taking $w = B^t$. Thus $dw
= B^t \log B\ dt$ or $\frac{dw}{w} =\log B\ dt$. As $e^{-2\pi i  tk}
= (B^{t/\log B})^{-2\pi i k} = w^{-2\pi ik /\log B}$ we have
\bea\label{eq:rwrwrwrewrw} F_B'(b) & \ = \ & \sum_{k=-\infty}^\infty
\int_0^\infty e^{-\beta w} \beta w \cdot w^{-2\pi i k / \log B}\
\frac{dw}{w} \nonumber\\ & = & \sum_{k=-\infty}^\infty
\beta^{2\pi i k/\log B} \int_0^\infty e^{-u} u^{-2\pi i k/\log B} du \nonumber\\
& = & \sum_{k=-\infty}^\infty \beta^{2\pi i k/\log B} \Gamma\left(1
- \frac{2\pi i k}{\log B}\right), \eea where we have used the
definition of the $\Gamma$-function: \be\label{eq:defGammafn}
\Gamma(s) \ = \ \int_0^\infty e^{-u} u^{s-1}\ du, \ \ \ {\rm Re}(s)
> 0. \ee As $\Gamma(1) = 1$ we have \be\label{eq:Fprimebgamma}
F_B'(b) \ = \ 1 + \sum_{m = 1}^\infty \left[ \beta^{2\pi im/\log B}
\Gamma\left(1 - \frac{2\pi i m}{\log B}\right) +\beta^{-2\pi im/\log
B} \Gamma\left(1 +\frac{2\pi i m}{\log B}\right)\right]. \ee

\begin{rek}\label{rek:beautyseriesexp}
The above series expansion is rapidly
convergent, and shows the deviations of $\log_B\zeta \bmod 1$ from
being equidistributed as an infinite sum of special values of a
standard function. As $\beta = B^b$ we have $\beta^{2\pi im/\log B}
= \cos(2\pi mb) + i \sin(2\pi mb)$, which gives a Fourier series
expansion for $F'(b)$ with coefficients arising from special values
of the $\Gamma$-function. \end{rek}

We can improve \eqref{eq:Fprimebgamma} by using additional
properties of the $\Gamma$-function. If $y\in \R$ then from
\eqref{eq:defGammafn} we have $\Gamma(1-iy) =
\overline{\Gamma(1+iy)}$ (where the bar denotes complex
conjugation). Thus the $m$\textsuperscript{th} summand in
\eqref{eq:Fprimebgamma} is the sum of a number and its complex
conjugate, which is simply twice the real part. We have
formulas for the absolute value of the $\Gamma$-function for large
argument. We use (see (8.332) on page 946 of \cite{GR}) that
\be\label{eq:gamma1+ix} |\Gamma(1+ix)|^2 \ = \ \frac{\pi
x}{\sinh(\pi x)} \ = \ \frac{2 \pi x}{e^{\pi x} - e^{-\pi x}}. \ee
Writing the summands in \eqref{eq:Fprimebgamma} as $2{\rm
Re}\left(e^{-2\pi i mb} \Gamma\left(1+\frac{2\pi i m}{\log
B}\right)\right)$, \eqref{eq:Fprimebgamma} becomes
\bea\label{eq:Fprimebgammanext} F_B'(b)  & = & 1 + 2\sum_{m =
1}^{M-1} {\rm Re}\left(e^{-2\pi i mb} \Gamma\left(1+\frac{2\pi i
m}{\log B}\right)\right)\nonumber\\ & & \ \ \ \ +\ 2\sum_{m =
M}^\infty {\rm Re}\left(e^{-2\pi i mb} \Gamma\left(1+\frac{2\pi i
m}{\log B} \right)\right). \eea The rest of the claims of Theorem
\ref{thm:main} follow from simple estimation, algebra and
trigonometry. \hfill $\Box$

With constants as in the theorem, if we take $M=1$ and $B = e$
(resp., $B=10$) the error is at most $.00499$ (resp., .378), while
if $M=2$ and $B=e$ (resp., $B=10$) the error is at most $3.16 \cdot
10^{-7}$ (resp., .006). Thus just \emph{one} term is enough to get
approximately five digits of accuracy base $e$, and two terms give
three digits of accuracy base $10$! For many bases we have reduced
the problem to evaluating ${\rm Re}\left(e^{-2\pi i b}
\Gamma\left(1+\frac{2\pi i}{\log B}\right)\right)$. This example illustrates the power of Poisson
Summation, taking a slowly convergent series expansion
and replacing it with a rapidly converging one.

\begin{cor}\label{cor:neverbenford} Let $\zeta$ have the standard exponential
distribution. There is no base $B>1$ such that $\zeta$ is Benford base
$B$.
\end{cor}

\begin{proof} Consider the infinite series expansion
in \eqref{eq:thmeqmainFprime}. As $e^{-2\pi i m b}$ is a sum of a
cosine and a sine term, \eqref{eq:thmeqmainFprime} gives a rapidly
convergent Fourier series expansion. If $\zeta$ were Benford base
$B$, then $F_B'(b)$ must be identically $1$; however,
$\Gamma\left(1+\frac{2\pi i m}{\log B}\right)$ is never zero for $m$
a positive integer because its modulus is non-zero (see
\eqref{eq:gamma1+ix}). As there is a unique rapidly convergent
Fourier series equal to $1$ (namely, $g(b) = 1$; see \cite{SS} for a
proof), our $F_B'(b)$ cannot identically equal 1.
\end{proof}


\section{Analyzing $N^{1-\delta}$ intervals simultaneously}\label{sec:analyzingintervalssimult}

We show why in addition to $\gep > 0$ we also needed $\gep > 1/3 - \delta/2$ when we analyzed $N^{1-\delta}$ intervals simultaneously in \eqref{eq:mknsimultepsilon}; we thank one of the referees for providing this detailed argument.

Let $Y_1, \dots, Y_N$ be iidrv with $\E[Y_i] = 0$, ${\rm Var}(Y_i) = \sigma^2$, $\E[|Y_i|^3] < \infty$, and set $S_N = (Y_1 + \cdots + Y_N) / \sqrt{N\sigma^2}$. Let $\Phi(x)$ denote the cumulative distribution function of the standard normal. Using a (non-uniform) sharpening of the Berry-Ess$\acute{{\rm e}}$en estimate (see, for example, \cite{Pe}), we find that for some constant $c > 0$ \be\label{eq:refeq1} \left|{\rm Prob}(S_N \le x) - \Phi(x)\right| \ \le \ \frac{c \E[|Y_1|^3]}{\sigma^3 \sqrt{N} (1+|x|)^3}, \ \ \ x \in \R, \ N \ge 1.\ee Taking $Y_i = w_i - N^{\delta-1}$, where $w_i$ is defined by
\eqref{eq:defnwi}, yields \bea S_N & \ = \ & \frac{M_N - N^\delta}{\sqrt{N^\delta (1 - N^{\delta-1})}} \nonumber\\ \sigma^2 &=& N^{\delta - 1} (1 - N^{\delta-1}) \nonumber\\ \E[|Y_i|^3] & \le & 2N^{\delta-1}. \eea Thus \eqref{eq:refeq1} becomes \be\label{eq:refeq2} \left|{\rm Prob}\left( \frac{M_N - N^{\delta}}{\sqrt{N^\delta (1-N^{\delta-1})}} \le x\right) - \Phi(x)\right| \ \le \ \frac{3cN^{-\delta/2}}{(1+|x|)^3}  \ee for all $N \ge N_0$ (for some $N_0$ sufficiently large, depending on $\delta$).

For each $N$, $k$ and $\gep$ consider the event \be A_{N,k,\gep} \ = \ \left\{ \frac{M_{k;N}-N^\delta}{\sqrt{N^\delta (1-N^{\delta-1})}} \ \in \ [-N^\gep, N^\gep]\right\}. \ee Then as $N\to\infty$ we have \be {\rm Prob}\left( \bigcap_{k=1}^{N^{1-\delta}} A_{N,k,\gep}\right) \ \to \ 1 \ee provided that \be\label{eq:refeqplus} \sum_{k=1}^{N^{1-\delta}} {\rm Prob}\left(A_{N,k,\gep}^{\rm c}\right) \ \to 0  \ee as $N\to\infty$. Using \eqref{eq:refeq2} gives \bea {\rm Prob}\left(A_{N,k,\gep}^{\rm c}\right) & \ \le \ & \frac{6cN^{-\delta/2}}{(1+N^\gep)^3} + 2\left(1 - \Phi(N^\gep)\right) \nonumber\\ & \le & 6cN^{-\delta/2 - 3\gep} + \sqrt{\frac{2}{\pi}} N^{-\gep} \exp(-N^{2\gep}/2) \eea (see, for example, \cite{Fe}). Thus the sum in \eqref{eq:refeqplus} is at most \be 6cN^{1-3\delta/2-3\gep} + \sqrt{\frac{2}{\pi}} N^{1-\delta-\gep} \exp(-N^{2\gep}/2), \ee and this is $o(1)$ provided that $\gep > 0$ \emph{and} $\gep > 1/3 - \delta/2$.


\ \\

\end{document}